\documentclass[12pt]{article}
\usepackage{amsfonts, latexsym, amsmath}

\font\fontA=cmr17

\font\mai=cmr12

\font\quadrat=lasy5 at 12pt
\def\quadr{\raise-1pt\hbox{\quadrat 2}}

\newtheorem{thm}{Theorem}
\newtheorem{rem}{Remark}
\newtheorem{lem}{Lemma}

\newcommand{\beq}{\begin{equation}}
\newcommand{\gb}{\goodbreak}
\newcommand{\eeq}{\end{equation}}

\newcommand{\al}{\alpha}
\newcommand{\be}{\beta}

\newcommand{\e}{\varepsilon}
\newcommand{\f}{\varphi}

\newcommand{\la}{\lambda}

\newcommand{\La}{\Lambda}

\newcommand{\Jn}{\langle \xi \rangle}
\newcommand{\JnA}{\langle \xi_1 \rangle}
\newcommand{\JnB}{\langle \xi_\nu \rangle}
\newcommand{\txi}{\tau(\xi)}

\newcommand{\gas}{\gamma^s}

\newcommand{\N}{\mathbb{N}}
\newcommand{\R}{\mathbb{R}}

\newcommand{\cinf}{{\mathcal C}^\infty}

\newcommand{\cm}{{\mathcal C}^m}
\newcommand{\A}{\mathcal A}
\newcommand{\ana}{\A}
\newcommand{\F}{\mathcal F}
\newcommand{\G}{\mathcal G}
\newcommand{\E}{\mathcal E}
\newcommand{\LL}{\mathcal L}
\newcommand{\Q}{\mathcal Q}
\newcommand{\CC}{\mathcal C}
\newcommand{\X}{\mathcal X}
\newcommand{\Qe}{Q_\e}

\newcommand{\Ee}{E_\e}

\newcommand{\wh}{\widehat}
\newcommand{\wt}{\widetilde}
\newcommand{\de}{\partial}
\newcommand{\dt}{\de_{t}}
\newcommand\dx{\de_{x}}

\newcommand{\ug}{\,=\,}
\newcommand{\mn}{\,\leq\,}
\newcommand{\mg}{\,\geq\,}
\newcommand{\+}{\,+\,}

\newcommand\noi{\noindent}
\newcommand\sm{\smallskip}
\newcommand\md{\medskip}
\newcommand\bg{\bigskip}
\newcommand\hb{\hbox}

\newcommand\edoc{\end{document}}
\newcommand\lb{\label}
\newcommand\rf{\ref}

\begin{document}
\title{\fontA Propagation of analyticity for a
class of nonlinear hyperbolic equations}
\author{\mai Sergio Spagnolo}
\date{}
\maketitle

\hfill {\it A Giovanni Prodi,} 

\hfill {\it indimenticabile scienziato, maestro ed amico}
\bigskip \bigskip 

\begin{abstract} \noi We consider the hyperbolic semilinear equations of the form 
$$
\dt^mu \+ a_1(t)\,\dt^{m-1}\dx u \+ \cdots \+ a_m(t)\,\dx^mu \ug f(u),
$$
$f(u)$ entire analytic, with characteristic roots  satisfying the condition 
$$
\la_i^2(t)+\la_j^2(t) \le M( \la_i(t)-\la_j(t))^2 , \ \quad \hb{for} \ \ i\neq j,
$$
and we prove that, if  the $a_h(t)$  are analytic, each solution  bounded in $\cinf$ enjoys the propagation of analyticity; while if $a_h(t)\in \,\cinf$, this property holds for those solutions which are bounded in some Gevrey class.
\end{abstract}

\section{Introduction} 

\noi The  linear operator
\beq\lb{N-sys}
\LL \,U\, =\, U_t\+ \sum_{h=1}^n A_h(t,x) \,U_{x_h}  
 \ \quad \hb{\rm on}\ \ [0,T]\times \R^n,
\eeq
where the $A_h$'s are $N\times N$ matrices, $U\in \R^N$,
 is {\it hyperbolic}
when, for all  $\xi \in  \R^n$,  
 the matrix 
$\,\sum A_h(t,x)\,\xi_h\,$ has  real eigenvalues $\la_ j(t,x,\xi)$,
$1\le j\le N$.

\noi Denoting by $\mu(\la)$ the multiplicity of the eigenvalue $\la$, we call {\it multiplicity} of  (\rf{N-sys}) the integer 
$
m = \max_{t,x,\xi}\,\max_j\{\mu (\la_j(t,x,\xi)\}.
$
The case $m=1$ corresponds to the {\it strictly hyperbolic systems}. \sm

We study the regularity of solutions to  nonlinear weakly hyperbolic system, in particular, {\it semilinear systems} 
\beq\lb{nonlin}
\LL\, U  \ug f(t,x,U)\,,
\eeq
where $\,U:[0,T]\times \R^n \to \R^N$, and $f(t,x,U)$ is a $\R^N$-valued,  analytic  function, typically a polynomial in the scalar components of $U$.

More precisely, assuming the coefficients of $\,\LL$ analytic  in $x$, we investigate under which additional assumptions a solution $U(t,x)$ of (\rf{nonlin}), analytic at the initial time, keeps its analyticity, i.e., satisfies 
\beq \lb{ana-pro}
U(0,\cdot) \in \ana(\R^n) \quad \Longrightarrow \quad U(t,\cdot) \in \ana (\R^n)
\quad \forall\, t\in [0,T]
\eeq
\gb

\noi Actualy, we consider two versions of  (\rf{ana-pro}), the first weaker and the second one stronger than (\rf{ana-pro}):
\begin{eqnarray}
\lb{L2-ana-pro} U(0,\cdot) \in \ana_{L^2}(\R^n) & \Longrightarrow & U(t,\cdot) \in \ana_{L^2}(\R^n) \quad \forall\, t\in [0,T]\,,
\\
\lb{loc-ana-pro} U(0,\cdot) \in \ana(\Gamma_0) \quad& \Longrightarrow & U(t,\cdot) \in \ana(\Gamma_t)  \ \qquad \forall\, t\in [0,T]\,,
\end{eqnarray}
where $\ana_{L^2}(\R^n)$ is the class of (analytic) functions $\f (x)\in H^\infty$  such that $\|\de^j\f\|_{L^2}\le C\La^j\,j!\,$, while $\Gamma$ is a {\it cone of determinacy} for the operator $\LL$ with base $\Gamma_0$ (at $t=0$) and sections $\{\Gamma_t\}$.
\sm

The propagation of analyticity is a natural property for nonlinear hyperbolic equations. Indeed, on one side, the theorem of Cauchy-Kovalewsky ensures the validity of  (\rf{ana-pro}) in some time interval $[0,\tau[$  (the problem is to prove that $\tau =T$), on the other side, by the Bony-Schapira's theorem, the Cauchy problem for any linear (weakly) hyperbolic system    is globally well posed in the class of   analytic functions.
\sm 

The first results of analytic propagation goes back to Lax (\cite{L},  1953) who considered (\rf{nonlin}) with $n=1\,$ in the strictly hyperbolic  case, and proved
(\rf{loc-ana-pro}) for those solutions  which are a priori bounded in $\,\CC^1$. Later on Alinhac and M\'etivier (\cite{AM}, 1984) extended this results to several space dimensions, but assuming that $U(t,\cdot)$ is bounded in $H^s(\R^n)$ for $s$ greater than some $\bar{s}(n)$.
\sm

In the weakly hyperbolic (nonlinear) case, the first results were  concerning
a second order equation of the form
\beq\lb{S}
\LL_0 \,u \, \equiv\, \sum_{i,j}^{1,n}\,\de_{x_i}(a_{ij}(t,x)\,\de_{x_j}u) \ug f(u),
\quad \sum a_{ij}\,\xi_i\xi_j \ge 0,
\eeq 
with $f(u), a_{ij}(t,x)$ analytic :
\sm

\noi {\bf Theorem A} \ (\cite{S}, 1989)

 {\it
\noi {\it i}) \ \ In the special case when $
a_{ij} = \be_0(t)\,\al_{ij}(x)$,  
a solution of $(\rf{S})$ enjoys $(\rf{loc-ana-pro})$ as long as remains bounded in $\cinf$.

\noi {\it ii}) \ \  In the general case, a solution  $u(t,\cdot)$ enjoys $(\rf{loc-ana-pro})$ provided it  is bounded in some Gevrey class $\gas$ with $s<2$.
} \sm

We recall that the Cauchy problem for any strictly hyperbolic linear system is globally  wellposed in $\cinf$. On the other hand, the Cauchy problem for the linear equation $\LL_0 u=0$, i is globally wellposed in $\cinf$ n the special case ({\it i}), whereas it is only globally wellposed in $\gas$ for $s<2$ in the general case ({\it ii}). 
Thus, it is natural to formulate the following
\sm

\noi {\bf Conjecture} \ {\it  In order to get the analytic propagation for a given solution  to a weakly hyperbolic system $\LL \,U=f(t,x,U)$, it is sufficient to assume {\rm a priori} that $\,U(t,\cdot)$ is bounded in some functional class $\mathcal X$ in which the Cauchy problem for the linear systems $\LL U+B(t,x)U=f(t,x)$  is globally well posed.}
 \sm
 
 \noi  [\,Typically the space $\X$ is equal to $\,\cinf$ or to some Gevrey class $\gas$ ]
\sm
 
 In the case when $\LL\,$ is a weakly hyperbolic operator of the {\it general type} (\rf{N-sys}), this Conjecture says that a solution $U(t,\cdot)$ enjoys the analytic propagation a long as  remains bounded in some Gevrey class $\gas$ of order $s<m/(m-1)$, where $m$ is the multiplicity of $\LL$. Indeed, Bronshtein's Theorem (\cite{B},  1979) states that, for any linear system $\LL\, U + B(t,x)U = f(t,x)$ with analytic coefficients in $x$, the Cauchy problem is  well-posed in these Gevrey classes. 
 \sm

Actually, this fact was proved in two special cases: time depending coefficients, and one space variable. More precisely:
\sm

\noi {\bf Theorem B} \  (\cite{DS}, 1999) 
\ {\it A solution of 
$$
U_t \+ \hb{$\sum_{j=1}^n A_j(t)$}\,U_{x_j} \ug f(t,x,U), \quad x\in\R^n,
$$
satisfies $(\rf{L2-ana-pro})$ as long as $\,U(t,\cdot)$ remains bounded in some $\gas$ with $s<m/(m-1)$. }
\sm

\noi {\bf Theorem C} \  (\cite{ST}, 2010) \ {\it A solution of 
$$
 U_t\+A(t,x)\,U_{x} \ug f(t,x,U), \qquad x\in \R,
$$
satisfies $(\rf{loc-ana-pro})$ as long as  $\,U(t,\cdot)$ remains bounded in some $\gas$ with $s<m/(m-1)$.} 
\md

\noi The study of the general case (coefficients depending on $(t,x)$, and $n\ge 2$)  is in progress.
\md

\noi {\bf Open Problem.} \  
To prove the sharpness of the bound $s<m/(m-1)$ in Theorems B or C. In particular: to construct a hyperbolic nonlinear system admitting a solution $U\in \cinf(\R^2)$ which is analytic on the halfplane $\{t<0\}$ but non analytic at some point of the line $t=0$. This kind of questions is related to the so called {\it Nonlinear Holmgren Theorem} (see \cite{M}). 
\md

\noi {\bf Acknowledgments.}\ We are indebted to Giovanni Taglialatela for his help to the drawing of this paper. 

\section{Main results}

\noi Hence, we consider  the scalar equations of the form
\beq\lb{eq-f}
\LL \,u \,\equiv\, \dt^mu \+ a_1(t)\,\dt^{m-1}\dx u \+ 
\cdots \+ a_m(t)\,\dx^mu \ug f(u) \,,
\eeq 
on $[0,T]\times \R$, where $f(u) = \sum_{\nu=0}^\infty u^\nu$ is an entire analytic, real function on $\R$.
We assume that the characteristic roots of the equation  are real functions, say
$$
\la_1(t) \le \la_2(t)\le  \ldots\le  \la_m(t)\,,
$$
which satisfying the condition
\beq \lb{diam}
\la_1^2(t)+\la_j^2(t) \mn M\,(\la_i(t)-\la_j(t)^2 , \ \  \quad \forall\, t\in [0,T]\quad (i\neq j).
\eeq

\begin{rem} Due to its symmetry  with respect to the roots $\la_j$, condition  $(\rf{diam})$ can be rewritten 
in term of the coefficients $\{a_h\}$ {\rm (Newton's theorem}. 
In particular {\rm (see \cite{KS})}: for a second order equation, $(\rf{diam})$ reads (for some $c>0$) 
$$
\Delta(t) \,\equiv\,  a_1^2(t)-4\,a_2(t) \mg c\, a_1^2(t) \,;
\lb{eq:C2}
$$
while for a  third order equation, it becomes
 \beq
\Delta(t) \mg c\, (a_1(t)a_2(t)-9\,a_3(t))^2 ,
\lb{eq:C3}
\eeq 
the discriminant being now
$ \Delta = 
   - 4\,a_2^3-27\,a_3^2+a_1^2\,a_2^2-4\,a_1^3\,a_3+18\,a_1a_2a_3$.
Particularly simple are the third order {\rm traceless} equations. 
i.e.,  when  $a_1 \equiv 0$:  here 
$a_2= -(\la_1^2+\la_2^2+\la_3^2)/2
\mn 0$, 
$\Delta = -4\,a_2^3\,-\,27\,a_3^2\,,$
so that  $(\rf{diam})$ becomes $\,\Delta \ge  -c\,a_2^{3}$, or equivalently  $\,\Delta \ge c\,a_3^2$. 
\end{rem}

Condition (\rf{diam}) for the linear equation $\LL u=0$ was introduced in [CO] as a sufficient (and almost necessary) condition for the wellposedness in $\cinf$.
A different proof of such a result, based on the quasi-symmetrizer, was given in \cite{KS}, where, also the case of non-analytic coefficients was considered:  it was proved that, if $a_h(t)\in \cinf([0,T])$ and (\rf{diam}) is fulfilled, then the Cauchy problem for 
$\LL u=0$ is well posed in each Gevrey class $\gas,\,s\ge 1$.

By these existence results, it is natural to expect some kind of  analytic propagation for the solutions 
which are bounded in $\cinf$ in case of analytic coefficients, or for those which are bounded in some Gevrey class $\gas$ in case of $\cinf$ coefficients.  

Actually, introducing the analytic, and Gevrey classes
\begin{eqnarray*}
 \ana_{L^2}  \!\!&=&\!\!
\big\{\f(x)\in\cinf(\R) : \|\de^j\f\|_{L^p(\R)} \mn C\,\La^j\,j!\,\big\}\,,
\\
 \gas_{L^2}  \!\! &=& \!\!
\big\{\f(x)\in\cinf(\R) : \|\de^j\f\|_{L^p(\R)} \mn C\,\La^j\,j!^s\big\}\,,
 \end{eqnarray*}
where $s\ge 1$, we prove:

\begin{thm}\lb{A} \
Assume that  the $a_j(t)$'s are analytic functions on $[0,T]$.
Then, for any solution of  $(\rf{eq-f})$ satisfying 
 \begin{eqnarray}\lb{PA1}
\sup_{0\le t\le T}\int_\R|\dt^{\,h}\dx^j u(t,x)|\,dx  \!\!&< &\!\!  \infty, \qquad \forall\, j\in\N,
\\
 \lb{PA2}
\dt^{\,h} u(0,\cdot) \!\!&\in &\!\! \ana_{L^2} ,
 \end{eqnarray}
 for $h=0, 1, \ldots, m-1$, it holds 
\beq \lb{PAA}
u \,\in \,\CC^{m-1}([0,T], \ana_{L^2}) \,.
\eeq
Under the same assumptions, we have also
\beq \lb{PAB}
 u\in \ana \,([0,T]\times \R) \,.
\eeq
\end{thm}

\begin{thm}\lb{B} \
If the $\,a_j(t)$'s are $\cinf$ functions on $[0,T]$,  the implication {\rm(\rf{PA2})} 
$\!\! \!\implies  \!\!\!$
{\rm(\rf{PAA})} holds true for those solutions which belong to 
$\,\cm\,([0,T], \gas_{L^2})$ for some $s\ge 1$\,.
\end{thm}

\noi {\bf Proof of Theorem 1.}
For the sake of simplicity, we shall perform the proof only in the case when the nonlinear term $f(u)$ is a monomial function, the general case requiring only minor additional computations. 
Thus, for a given integer $\nu\ge 1$, we consider the equation
\beq\lb{eq}
\dt^mu \+ a_1(t)\,\dt^{m-1}\dx u \+ 
\cdots \+ a_m(t)\,\dx^mu \ug u^\nu . 
\eeq 

\noi Putting
$$
\wh u(t,\xi) \ug \int_{-\infty}^{+\infty}e^{-i\xi x}\, u(t,x)\, dx,
$$
\beq\lb{V}
V(t,\xi) =\begin{pmatrix}
(i\xi)^{m-1}\,\wh u \cr (i\xi)^{m-2}\,\wh u\,' 
\cr \vdots \cr \wh u^{(m-1)}
 \end{pmatrix}\,, \quad 
  F(t,\xi) = \begin{pmatrix}0 \cr 0 \cr \vdots \cr f(t,\xi)
 \end{pmatrix} \,,
\eeq
and
\beq\lb{A(t)}
A(t) = \begin{pmatrix} \!0 & 1 & & 
\cr
&  \ddots &\ddots &
\cr
&&0&\!\!1
\cr
\,a_m(t) & \cdots&a_2(t) &\,a_1(t)
 \end{pmatrix} ,
\eeq
we transform equation (\ref{eq}) into the ODE's system
\beq \lb{sys}
V'+\, i\,\xi A(t)V \ug F(t,\xi),
\eeq
where
\beq \lb{F}
 f(t,\xi) \ug \underbrace{\wh u * \cdots * \wh u \,}_{\nu}.
\eeq 
\gb

\noi Our target is to prove that, if 
\begin{eqnarray}
\lb{FAP1}
\int_\R|\xi|^j\,|V(t,\xi)|\, d\xi \!\!&\le& \!\!  K_j < \infty
\ \quad \forall\, j\,, \quad \forall\,t\in [0,T],
\\
\lb{FAP2}
\int_\R|\xi|^j\,|V(0,\xi)|\, d\xi \!\!&\le& \!\!  C\,\La^jj!  \  \qquad \forall\, j,
\end{eqnarray}
then, for some new constants $\wt C,\, \wt \La$, it holds 
\beq\lb{FAPA}
\int_\R|\xi|^j\,|V(t,\xi)|\,d\xi\mn   \wt C\,\wt \La^jj! \,, \qquad \forall\, j
\,, \quad \forall\,t\in [0,T].
\eeq
Indeed, (\rf{FAP2})  is an easy consequences of (\rf{PA2}); while (\rf{PA1}) implies that $\{\dt^h \dx^j u(t,\cdot)\}$ is bounded in
$L^{\infty}(\R)$ for all $j$,  whence
  (\rf{FAP1}). 
     Finally,  taking into account that  $|V(t,\xi)|\le K <\infty$ (by (\rf{PA1})), we see that (\rf{FAPA}) implies (\rf{PAA}). 
     \sm

\noi To get this target, we firstly prove an apriori estimate for the {\it linear system} (\rf {sys}), without taking (\rf{F}) into account.
We follow \cite{KS}, but some modifications are needed in order to get an estimate suitable to the nonlinear case. The main tool   is  the theory of quasi-symmetrizer developed in \cite{J} and 
\cite{DS}.  
\sm

\noi {\bf Recalls on quasi-symmetrizer.} 
\sm

\noi{\bf \cite{DS}} :  For any matrix of the form (\rf{A(t)}) with real eigenvalues, we can find a family of Hermitian matrices  
\beq \lb{Qeps} 
\Qe(t) \ug \Q_0(t)+\e^2\Q_1(t)+\cdots +\e^{2(m-1)}\Q_{m-1}(t)
\eeq
such that the entries of the $\Q_r(t)$'s are polynomial functions of the coefficients $a_1(t), \ldots, a_m(t)$ (in particular inherit their regularity in $t$), and 
\begin{eqnarray}\lb{QS1}
& C^{-1}\e^{2(m-1)}\,|V|^2  \mn (\Qe(t)V, V) \mn  C\,|V|^2&
\\
\lb{QS2} & \big |\big(\Qe(t)A(t)-A(t)\Qe(t))V,\,V \big) \big | \mn C\,\e^{1-m}\,(\Qe(t)V,V).&
\end{eqnarray}
for all $\,V\in \R^m, \, 0<\e\le 1$.
\sm

\noi 
{\bf \cite{KS}} : If the eigenvalues of $A(t)$ satisfy the condition {\rm (\rf{diam})}, then $\Qe(t)$ is a {\it nearly diagonal matrix}, i.e., it satisfies,  for some constant $c>0$, independent on $\e$, \beq \lb{ND}
\qquad(\Qe(t)V,V) \mg c\, \sum_{j=1}^m q_{\e, jj}(t)v_j^2\,,  \qquad \forall\, V\in \R^m,
\eeq
 where $q_{\e, ij}$ are the entries of $\Qe$, $v_j$ the scalar components of $V$. \hfill$\Box$
\sm

In our assumptions,  the $a_h(t)$'s are analytic functions on $[0,T]$, consequently also
the entries  $q_{r,ij}(t), \, 1 \le i,j \le m$ of  the matrix $\Q_r(t)$ will be analytic.
Therefore,
putting together all the isolated zeroes of these functions, we form a partition of $[0,T]$,
 independent on $\e$,
\beq \lb{part}
0=t_0<t_1< \cdots <t_{N-1} < t_N=T,
\eeq
 such that, for each  $r, i, j$, it holds:
$$
\hbox{\rm either}\ \ q_{r,ij} \equiv 0, \quad \hbox{\rm or} \ \   \ q_{r,ij}(t) \neq 0 \ \  \ \forall t \in 
I_h = [t_{h-1}, t_h[\,.  
$$
\sm
 Now,  let us notice that, by  Cauchy-Kovalewsky, if at some point $t$ a solution to (\rf{eq}) satisfies $\dt^{\,h} u(t,\cdot) \in \ana_{L^2}(\R)$ for all $\,h \le m-1$,  then the same holds in a right neighborhood of $t$.
 Thus,  it will be sufficient to put ourselves inside one of the intervals $I_1, \ldots, I_N$. In other words it is not restrictive to assume that, for each $r,i,j$, 
\beq \lb{either}
\hbox{\rm either}\ \ q_{r,ij} \equiv 0, \quad \hbox{\rm or} \ \   \ q_{r,ij}(t) \neq 0 
\ \quad \hb{\rm for} \ \  0\le t <T.  
\eeq
Therefore, by  the analyticity of  $q_{r,ij}(t)$  we easily derive that
\beq \lb{q'}
|q_{r,ij}'(t)| 
\mn 
\frac C {T-t} \ |q_{r,ij}(t)|  \  \quad \hb{on} \ \  [0,T[.
\eeq

\noi  Next, following \cite{KS}, for any fixed $\xi\in\R$ we prove two different apriori estimates for a solution $V(t,\xi)$ of (\rf{sys}):  
 a {\it Kovalewskian} estimate in a (small) left neighborhood of $T$, $[T-\tau, \,T[$, and a {\it hyperbolic} estimate on $[0,\tau]$.

\noi [\,In the following  $C, C_j$ will be constants depending on the coefficients 
of  (\ref{eq})\,]
\begin{lem}
 Let $V(t,\xi)$ be a solution  of  $(\rf{sys})$ on  $[0,T[$, and put 
 \beq \lb{E-eps}
\Ee(t,\xi) \ug (\Qe(t)\,V(t,\xi), V(t,\xi)).
\eeq
Then,
 for any fixed $\xi\in\R$, the following estimates hold:
 \beq \lb{kov-est}
\dt\,|V(t,\xi)| \mn \frac{C_0} T\, |\xi| \,|V(t,\xi)| + |F(t,\xi)|, 
\eeq
\beq \lb{hyp-est}
\dt\sqrt {\Ee(t,\xi)} \mn  C_0\,\Big(\,\frac 1{T-t}+\e\, |\xi|\Big) \,\sqrt{\Ee(t,\xi)}\+  C_0\,|F(t,\xi)| ,  
\eeq
$C_0$ a constant depending only on  the coefficients of the equation, and on $T$.\\
\noi In particular, putting  \beq\lb{E-star}
\qquad \qquad E_* \ug E_{\e_*}\,,  \qquad \hb{\rm where} \ \  \e_* \ug \Jn^{-1}, \ \ \Jn \ug 1+|\xi|\,,
\eeq
$(\rf {hyp-est})$ gives 
\beq \lb{hyp-est1}
(\sqrt E_*)' \mn  C_0\,\Big(\,\frac 1{T-t}+1\Big) \,\sqrt E_*\+  C_0\,|F(t,\xi)|. 
\eeq
\end{lem}

\noi {\bf Proof:}
As an easy consequence of (\rf{sys}), we  get (\rf{kov-est}) with
$$
C_0\mg \max_{t\in [0,T]}\|A(t)\| \,, \ \  C_0 \mg 1\,.
$$  
To prove (\rf{hyp-est}) we differentiate (\rf{E-eps}) in time. 
Recalling (\rf{QS1}) we find 
\begin{eqnarray*}
\Ee'(t,\xi) & \ug& (\Qe'V,V) + (\Qe V',V) + (\Qe V,V') \\
& \ug& 
(\Qe'V,V) + i\,\xi\,((\Qe A - A^*\Qe)V,V)+ 2\,\Re(\Qe F,V) \\
&\le & 
K_\e (t,\xi)\,\Ee(t,\xi) +  C_1\,|F(t,\xi)|\,\sqrt{\Ee(t,\xi)}  
\end{eqnarray*}
where $V=V(t,\xi)$ and 
\beq \lb{Keps}
K_\e (t,\xi) \ug  \frac{|(\Qe' V,V)|}{(\Qe V,V)}\+ |\xi|\,\frac{|((\Qe A - A^*\Qe)V,V)|}
{(\Qe V,V)}. 
\eeq

\noi We  have to prove that
\beq \lb{K-est}
K_\e(t,\xi) \mn C\,\Big(\,\frac 1{T-t}+\e\, |\xi|\Big) \qquad\forall\, t\in [0, T[ \,. 
\eeq
Let us firstly note that the second quotient in  (\rf {Keps}) is estimated by $C\e$ by  the property (\rf{QS2})  of our quasi-symetrizer.   
To estimate the first quotient, apply to the nearly diagonality of the matrix $\Qe(t)$, i.e., (\rf{ND}): recalling (\rf{Qeps}),  and
noting that 
$|q_{r,ij}| \le \sqrt{q_{r,ii}\,q_{r,jj}}$ (since $\Q_r(t)$ is a symmetric matrix $\ge 0$), it follows
\begin{eqnarray*}
|(\Qe'V,V)| \!\!& \le & \!\! \sum_{r=0}^{m-1}± \e^{2r}\,\sum_{ij}^{1,n}|q_{r,ij}'||v_iv_j| \mn C\,(T-t)^{-1}
\sum_r \e^{2r}\,\sum_{ij}|q_{r,ij}||v_iv_j|  
\\
 \!\!& \le & \!\!  C\,(T-t)^{-1} \sum_r \e^{2r}\,\sum_{j}q_{r,jj}\,v_j^2 
 \ug C\,(T-t)^{-1} q_{\e, jj}\, v_j^2
 \\
 \!\!& \le & \!\! C_1\,(T-t)^{-1}(\Qe V,V)\,.
\end{eqnarray*}
This completes the proof of (\rf{K-est}), hence of (\rf{hyp-est}).
\hfill  $\Box$
\md

Next,   we define 
\begin{eqnarray}
\lb{tau}
& \tau(\xi) \ug  T - {|\xi|}^{-1}\,,\qquad \qquad \qquad \qquad \ \ &
\\
\lb{Phi}
&\Phi(t,\xi) =   C_0\,\min \big\{(T-t)^{-1}\!+1, \Jn \big\}
=
 \begin{cases}{\displaystyle C_0\,\big\{ (T-t)^{-1}+1\big\} \  \ 
\hb{\rm on} \ [\,0, \tau(\xi)]}
 \\
 \\
 {\displaystyle {C_0}\,\Jn \qquad \hb{\rm on} \   [\txi, T\,[ }
\end{cases}&
\\
\lb{rho}
&\rho(t,\xi) \ug  \int_t^T \Phi(s,\xi)\,ds \,.\qquad \qquad \qquad \qquad&
\end{eqnarray}
Therefore,  by (\rf{kov-est}) and (\rf{hyp-est1}) it follows
\begin{eqnarray} \nonumber
\big\{|V(t,\xi)|\big\}' \!\! &\le& \!\! \Phi(t,\xi)\, |V(t,\xi)| \ \+ C_0\, |F(t,\xi)| \quad \  \,
\hb{\rm on} \ \ [\txi, T[ \\
\lb{stimina} \big\{\sqrt {E_*(t,\xi)}\big\}' \!\! &\le&\!\!  \Phi(t,\xi) \,\sqrt {E_*(t,\xi)} \+ C_0\, |F(t,\xi)|   \quad \  
\hb{\rm on} \ \, [\,0,\txi\,[  
\end{eqnarray}
and thus, since $\,\rho'=-\Phi$, 
\begin{eqnarray*} 
\de_t\Big\{e^{\rho(t,\xi)}\,|V(t,\xi)| \Big\}
\!\! &\le& \!\!   C_0\,e^{\rho(t,\xi)}\,|F(t,\xi)| \quad \  \hb{\rm for} \ \ \txi \le t \le  T
\\
\de_t\Big\{e^{\rho(t,\xi)}\sqrt{E_*(t,\xi)}\Big\}
\!\! &\le& \!\!    C_0\,e^{\rho(t,\xi)}\,|F(t,\xi)| \quad \  \hb{\rm for} \ \ 0 \le t\le \txi \,.
\end{eqnarray*}

\noi By integrating in time, we find (omitting $\xi$ everywhere)
\begin{eqnarray}
\lb {est-V}
e^{\rho(t)}\,|V(t)| \!\! &\le& \!\!  e^{\rho(\tau)}\, |V(\tau)| \,+ C_0\int_{\tau} ^{\, t} 
e^{\rho(s)} |F(s)|\,ds
\\
\lb{est-E}
e^{\rho(\tau)}\sqrt{E_*(\tau)} \!\! &\le& \!\!  e^{\rho(0)} \sqrt{E_*(0)} \,+
C_0 \int_{0} ^{\tau} 
e^{\rho(s)} |F(s)|\,ds
\end{eqnarray}

\noi Now, by (\rf{QS1}) with $\e=\Jn^{-1}$, we know that
$$
C^{-1}\Jn^{-2(1-m)}\,|V(t,\xi)|^2 \mn {E_*(t,\xi)} \mn C\,|V(t,\xi)|^2,
$$ 
hence  we derive, form (\rf{est-V}) and (\rf{est-E}),
\begin{eqnarray*}
e^{\rho(t)}\,|V(t)| \!\!&=&\!\! C_1 \Jn^{m-1}\,e^{\rho(\tau)}\sqrt{E_*(\tau)} \+ C_0\int_{\tau} ^{\, t} 
e^{\rho(s)} |F(s)|\,ds \qquad \qquad \qquad
\\
\!\!&\le&\!\! C_1 \Jn^{m-1}\,\Big\{e^{\rho(0)} \sqrt{E_*(0)} + \int_{0} ^{\tau} 
e^{\rho(s)} |F(s)|\,ds \Big\}
+ C_0\int_\tau^t e^{\rho(s)} |F(s)|\,ds
\\
\\
\!\!&\le&\!\! C_2\, \Jn^{m-1}\,\Big\{e^{\rho(0)}\sqrt{E_*(0)} \,+
C_0 \int_{0} ^{t} 
e^{\rho(s)} |F(s)|\,ds\Big\}.
\end{eqnarray*}

\noi Recalling the definitions (\rf{Phi}) and  (\rf{rho}) of $\Phi$ and $\rho$, we get

$$
\rho(0,\xi ) \ug \int_0^T\! \Phi(s,\xi)\,ds 
\mn C_0\int_0^{\tau(\xi)}\!\! \Big\{\frac 1 {T-t} +1\Big\}\,dt \+ (T-\tau(\xi))\Jn
$$
and hence we derive, since $\,\dt \rho <0\,$ and $\,\tau(\xi)=T-|\xi|^{-1}$,  
 \beq
 \lb{est-rho}
\rho(t, \xi) \mn C\,\left(\log \Jn +1\right) \qquad \hb{\rm for all} \ \ t\in [0,T].
\eeq
Therefore we obtain, for some integer $N$, 
\beq\lb{stima}
e^{\rho(t, \xi)}|V(t,\xi)| \mn C\,\Jn^N\,|V(0,\xi)| \+ C\, \Jn^{m-1}\, \int_{0} ^{t} 
e^{\rho(s, \xi)} |F(s,\xi)|\,ds.
\eeq
By the way, we note that the last inequality ensures the wellposedness in $\cinf$ of the Cauchy problem  for the linear system (\rf{sys}).
\sm

Let us go back to the nonlinear equation $\,\LL u=u^\nu$. For our purpose we must consider a more general equation, namely
$$
\LL u \ug u_1\cdots u_\nu,
$$
where the $u_j=u_j(t,x)$ are given functions (actually, some $x$-derivatives of $u$).\\
In such a case, the function $F$ in (\rf{sys}) is 
\beq \lb{F-conv}
F(t,\xi) \ug \wh u_1* \cdots * \wh u_\nu ,
\eeq
where the convolutions are effected w.r. to $\xi$, and thus 
$$
|F(t,\xi)| \mn  \int_{\xi_1+\cdots +\xi_\nu=\xi}|\wh u_1(t,\xi_1) \cdots
 \wh u_\nu(t,\xi_\nu)|\, d\sigma_{(\xi_1, \ldots, \xi_\nu)}\, .
$$
We notice that the function $\, \xi \mapsto\min\{C, |\xi|\}$ is a sub-additive; consequently
for each fixed $t$ (see (\rf{Phi}),(\rf{rho})) the function $\Phi(t,\xi)$, hence also $\rho(t,\xi)$, is sub-additive in $\xi$.
On the other hand,  $\xi \to\Jn$  is sub-multiplicative. \\
Thus one has, for $\xi=\xi_1+\cdots+\xi_\nu$, 
\begin{eqnarray*}
& \rho(t,\xi) \mn\rho(t,\xi_1) + \cdots + \rho(t,\xi_\nu), \qquad
\Jn ^{m-1} \le \langle{\xi_1}\rangle^{m-1} \cdots \langle{\xi_\nu}\rangle ^{m-1},&
\\
&e^{\rho(t,\xi)}\,\Jn ^{m-1}\mn e^{\rho(t,\xi_1)}\,\JnA ^{m-1}\cdots\, e^{\rho(t,\xi_\nu)}\,\JnB ^{m-1},&
\end{eqnarray*}
whence, by (\rf{F-conv}), it follows the pointwise estimate
$$
e^{\rho}\,\Jn ^{m-1}|F|  \mn \big (e^{\rho}\,\Jn ^{m-1}|\wh u_1| \big) * \cdots * \big(e^{\rho}\,\Jn ^{m-1}|\wh u_\nu| \big) .  
$$
 Now, if $\,V_j(t,\xi)$ are the vectors formed as  $V(t,\xi)$ (see (\rf{V})), with $u_j$ in place of $u$, we have
 $$
\Jn^{m-1}\,|\wh u_j(t,\xi)| \mn |V_j(t,\xi)|\,, \qquad j=1,\ldots,\nu,
$$ 
and thus, going back to (\rf{stima}), we obtain
$$
e^{\rho(t, \xi)}|V(t,\xi)| \mn C\Jn^N\,|V(0,\xi)| \+ C\int_{0} ^{t} 
\big(e^{\rho}\,|V_1|  * \cdots * e^{\rho}\,|V_\nu| \big)(s,\xi)\,ds.
$$
 Finally,  we integrate in $\xi\in \R$ to get
\beq \lb{NLest}
\E(t,u) \mn C\int_\R |V(0,\xi)| \,\Jn^Nd\xi \+ C \int_0^t \E(s,u_1) \cdots \E(s,u_\nu)  \,ds\,,
\eeq
where we define the $\,\cinf$-{\it energy}
\beq \lb{energy}
\E (t,u) \ug \int_\R e^{\rho(t, \xi)}|V(t,\xi)|\, d\xi.
\eeq
We emphasize that, by virtue of our assumption $(\rf{FAP1})$,   and $(\rf{est-rho})$, we have
\beq\lb{M-0} \E(t,u) \mn 
M_0  \, <\, \infty \qquad (0\le t\le T).
\eeq

\noi Differentiating $j$ times in  $x$ the equation $\LL u=u^\nu$, we get
$$
\LL(\de^j u)\ug j! \sum_{h_1+\cdots+h_\nu=j} \frac{\de^{h_1}u \cdots \de^{h_\nu}u}
{h_1! \cdots h_\nu!} \qquad (\hb{\rm where} \ \ \de = \dx),
$$
and to this equation we apply the estimate (\rf{NLest}) with $\,u_j=\de^ju$. We obtain:
\beq \lb{STIMA}
\frac{\E_j(t)}{j!} \mn C\int_\R \frac {|V_j(0,\xi)|}{j!}\,\Jn^N d\xi \+ C \sum_{|h|=j}\int_0^t  
\frac{\E_{h_1}(s)}{h_1!}\cdots \frac{\E_{h_\nu}(s)}{h_\nu!} \,ds,
\eeq
where $V_j(t,\xi)$ is the vector associated to $u_j \equiv \de^ju$, and 
$$
\E_j(t) \ug \E(t, \de^ju). 
$$
Putting
$$
\alpha_j(t) \ug \int_\R {|V_j(0,\xi)|}\,\Jn^N d\xi 
\+  j!\, \sum_{|h|=j}\int_0^t  
\frac{\E_{h_1}(s)}{h_1!}\cdots \frac{\E_{h_\nu}(s)}{h_\nu!} \,ds \,,
$$
we rewrite (\rf{STIMA}) as
\beq\lb{STIMA-1}
\E_j(t) \mn C\, \alpha_j(t) \,.
\eeq
Next, we  introduce the {\it super-energies}
\beq \lb{F-cal}
\F(t,u) \ug \sum_0^\infty\,{\E_j(t)}\, \frac{r(t)^j}{j!}\,,
\eeq
\beq\lb{G-cal} 
\G(t,u) \ug  \sum_0^\infty\, \alpha_j(t)\,\frac{r(t)^j}{j!}\  ,  
\quad
\G^1(t,u) \ug  \sum_1^\infty\,{\alpha_j(t)}\, \frac {r(t)^{j-1}}{(j-1)!}\  ,
\eeq
where $r(t)$ is a decreasing, positive function on [0,T] to be defined
later.

\noi By differentiating in time, we find 
\begin{eqnarray*}
 \G'  \!\!&=&\!\! \sum_0^\infty \,\alpha_j'\, \frac{r^j}{j!} + \sum_1^\infty\, \alpha_j\, \frac{r^{j-1}}{(j-1)!}\, r'
\ug \sum_{j=0}^\infty\sum_{|h|=j}\,{\E_{h_1}}\frac {r^{h_1}}{h_1!}\cdots{\E_{h_\nu}}
\frac {r^{h_\nu}}{h_\nu!}+ r'\,\G^1 \\
\!\!&=&\!\!  \Big\{\sum_{h=0}^\infty {\E_{h}}\frac {r^{h}}{h!}\Big\}^\nu+
r'\,\G^1 \ug \F^\nu \+ r'\,\G^1 ,
\end{eqnarray*}
and hence, noting that $\F(t)\le C\,\G(t)$ by (\rf{STIMA-1}), 
\beq \lb{est-G'}
\G' \mn C^\nu \,\G^\nu \+ r'\,\G^1.
\eeq
Now, noting that (by (\rf{FAP1}) and (\rf{M-0}))
$$
\alpha_0(t)  \ug \int_R|V(0,\xi)|\,\Jn^N d\xi \,+ \int_0^t\E(s)\,ds \mn K_N+M_0  \,\equiv\, M,
$$
by  the definition (\rf{G-cal}) of $\G(t)$ it follows
$$
\G(t) \mn \alpha_0(t) \+ r(t)\,\G^1(t) \mn M \+ r(t) \, \G^1(t)\,.
$$
From this inequality it follows, arguing by induction w.r. to $\nu$,
$$
\G^\nu  \mn M^\nu \+  r\, \G^1\,\big(\G+M\big)^{\nu-1};
$$
consequently (\rf{est-G'}) gives (putting $\phi(\G)=C^\nu (M + \G)^{\nu-1}$)  
\beq \lb{est-G-top}
\G' \mn \G^1\big\{r'+ r\,\phi(\G)\big\}\+(CM)^\nu.
\eeq
On the other hand, by virtue of our assumption (\rf{FAP2}),  we see that
$$
\G(0,u) \ug \sum_{j=0}^\infty  \Big\{\int_\R {|V_j(0,\xi)|}\,\Jn^N\, d\xi\Big\}\,
\frac{r(0)^j}{j!} \, <\, \infty\,. 
$$
provided $\,r(0)\equiv r_0$ is small enough.
Therefore, taking
\beq\lb{L,r}
L\ug \G(0,u)+(CM)^\nu\, T \, , \qquad r(t) \ug r_0\,e^{-\phi(L)\,t}\,, 
\eeq
we can derive from (\rf{est-G-top}) the estimate
\beq\lb{G<L}
\G(t,u) \,<\, L \qquad \hb{\rm for all} \ \ t\in [0,T].
\eeq

\noi {\bf  Proof of (\rf{G<L}). }  Since $L>\G(0)$, this estimate holds true in a right neighborhood of $t=0\,$ by  Cauchy-Kovalewsky. Then  assuming  that, for some $\tau_*<T$, (\rf{G<L}) holds for all
$\, t <\tau_*$ but not at $t=\tau_*$,  we have $\G(\tau_*)=L$, and hence, with $r(t)$ as in (\rf{L,r}), 
$$
r'(t) +r(t)\,\phi(\G(t)) \mn r'(t)+r(t)\, \phi(L) \mn 0 \qquad \hb{\rm on} \ \ [0,\tau_*[ \,.
$$
This yelds a contradiction; indeed, by (\rf{est-G-top}), $$
\G(t) \mn \G(0)+(CM)^\nu \tau_* \,<\, L \qquad \hb{\rm on} \ \ [0,\tau_*] \,.
$$
{\bf Conclusion of the Proof of Theorem 1.} 
Recalling that $\F(t,u) \le C\,\G(t,u)$, (\rf{G<L}) says  that 
$\F(t,u)<CL$ on $[0,T]$. Therefore, by (\rf{F-cal}),  we get our goal  (\rf{FAPA}):   
\begin{eqnarray*}
\int_\R|V(t,\xi)|\,|\xi|^j\,d\xi 
&\!\!\le\!\!&
\int_\R e^{\rho(t,\xi)}|V(t,\xi)|\,|\xi|^jd\xi
\ug \E(t,\de^ju) \mn \F(t)\,r(t)^{-j}\,j! 
\\
&\!\!\le\!\!& CL\, \big\{r_0\,e^{\phi(L)T}\big\}^j \,j!
\ug \wt C\,\wt \Lambda^{\,j+1}\,j!\,.
\end{eqnarray*} 

To prove (\rf{PAB}), i.e.,  the global analyticity of the solution $u$ in $(t,x)$, it is sufficient to resort to Cauchy-Kovalewski.
\hfill $\Box$

\md
\begin{rem} \lb{r_eta} \ The previos proof of $(\rf{G<L})$ is somewhat formal, since it  assumes  not only that $\G(t)<\infty$, but also that $\G^1(t)< \infty $ on $[0,\tau_*[$. To make the proof more precise we must replace the radius function $r(t)$ by $r_\eta(t)=\eta\exp(-\phi(L)t)$, $\eta<1$, and apply the previous computation to the corresponding functions $\,\G_\eta(t)$ and $\,\G^1_\eta(t)$. Finally let $\eta \to 1$ {\rm (see [ST] for the details)}.   
\end{rem}
\md

\noi {\bf Proof of Theorem 2.} \
 The proof is not very different from that of Thm.1, thus we give only a sketch of it.

The main difference is that the entries $q_{r,ij}(t)$ are no more analytic, but only $\cinf$, hence (\rf{q'})
fails. 
However, for any function $f\in \CC^k([0,T])$ it holds
$$
|f'(t)| \mn \La(t)\, |f(t)|^{1-1/k}\, \|f\|_{\CC^k([0,T])},
$$
for some $\La\in L^1(0,T)$ [this was proved in \cite {CJS} in the case $f(t)\ge 0$, and in  \cite{T} in the general case].
Therefore, recalling  that $\Qe(t)$ is a nearly diagonal matrix, and proceeeding
 as in \cite{KS}, we get, for all integer $k\ge 1$, 
\beq\lb{Q'}
 |\big (\Qe'(t)V(t,\xi),V(t,\xi) \big)| \mn \Lambda_k (t) \,\big (\Qe(t)V(t,\xi),V(t,\xi) \big)^{1-1/k}\, |V(t,\xi)|^{2/k}
 \eeq
for some $\Lambda_k\in L^1(0,T)$, independent of $\e$. 
Differently from Thm. 1, we need now to consider only the {\it hyperbolic energy}
$$
E_*(t,\xi) \ug (Q_{\e^*}(t)V,V) \quad \hb {\rm with}\  \  \e = |\xi|^{-1}.
$$
Thanks to (\rf{Q'}),  we prove  (for every integer $k\ge 1$) the estimate
$$
\big\{\sqrt {E_*(t,\xi)}\big\}' \mn C_0\, \Phi(t,\xi) \,\sqrt {E_*(t,\xi)} \+ C_0\, |F(t,\xi)| 
$$
 on all the interval $[0,T]$, where 
$$
\Phi(t,\xi) \ug\La_k(t) |\xi|^{2(m-1)/k} \+ 1
$$
Note that $\Phi$ is sub-additive w.r. to $\xi$ as soon as $k\ge 2(m-1)$.

\noi Next, putting
$$
\rho(t,\xi) \ug \int_t^T \Phi(t,\xi)\, d\xi  \,\equiv\, |\xi|^{2(m-1)/k}\int _t^T \La_k(s)\,ds \+  (T-t) \,,
$$
we define the {\it Gevrey-energy}
$$
\E(t,u) \ug \int _\R\, e^{\rho(t,\xi)}\sqrt {E_*(t,\xi)}\,d\xi .
$$
We conclude as in the proof of Thm.1.

\bg\bg

\noi {\small S. Spagnolo \\
   Dept. of  Mathematics L.Tonelli, University of Pisa \\
    Largo B.~Pontecorvo~5, 56127 Pisa, Italy \\
    {\it spagnolo@dm.unipi.it}}

\edoc